\documentclass[12pt]{article}
\usepackage{amsmath}
\usepackage{amsfonts, amssymb}

\begin{document}
\bibliographystyle{plain}
\pagenumbering{arabic}
\raggedbottom

\newtheorem{theorem}{Theorem}[section]
\newtheorem{lemma}[theorem]{Lemma}
\newtheorem{proposition}[theorem]{Proposition}
\newtheorem{corollary}[theorem]{Corollary}
\newtheorem{conjecture}[theorem]{Conjecture}
\newtheorem{definition}[theorem]{Definition}
\newtheorem{example}[theorem]{Example}
\newtheorem{condition}{Condition}
\newtheorem{main}{Theorem}

\setlength{\parskip}{\parsep}
\setlength{\parindent}{0pt}

\hfuzz12pt

\newcommand{\cl}[1]{{\mathcal{C}}_{#1}}
\newcommand{\corr}{\mbox{Corr}}
\newcommand{\cov}{\mbox{Cov}}
\newcommand{\der}[3]{\ra{#1}{#2}_{#3}}
\newcommand{\expo}{\mbox{Exp}}
\newcommand{\tr}[1]{\mbox{\rm{tr}}\left(#1 \right)}
\newcommand{\var}{\mbox{Var}\;}
\newcommand{\supp}{\mbox{supp}}
\newcommand{\tends}{\rightarrow \infty}
\newcommand{\ep}{{\Bbb {E}}}
\newcommand{\pr}{{\Bbb {P}}}
\newcommand{\re}{{\mathcal{R}}}
\newcommand{\rt}{\widetilde{\rho}}
\newcommand{\vt}{\rat{\widetilde{V}}}
\newcommand{\bds}{\begin{displaystyle}}
\newcommand{\eds}{\end{displaystyle}}
\newcommand{\vc}[1]{{\mathbf{#1}}}
\newcommand{\ra}[2]{{#1}^{(#2)}}
\newcommand{\rat}[1]{\ra{#1}{\tau}}
\newcommand{\vd}[1]{{\boldsymbol{#1}}}
\newcommand{\sif}[2]{{\cal{M}}_{#1}^{#2}}
\newcommand{\fish}{J_{\rm st}}

\def \outlineby #1#2#3{\vbox{\hrule\hbox{\vrule\kern #1% 
\vbox{\kern #2 #3\kern #2}\kern #1\vrule}\hrule}}%
\def \endbox {\outlineby{4pt}{4pt}{}}%

\newenvironment{proof}
{\noindent{\bf Proof\ }}{{\hfill \endbox
}\par\vskip2\parsep}

\newlength{\dsl}
\newcommand{\doublesub}[2]{\settowidth{\dsl}{$\scriptstyle
#2$}\parbox{\dsl}{\scriptsize\centering { \normalsize $\scriptstyle #1$}
\\ {\normalsize $\scriptstyle #2$}}}

\title{An information-theoretic Central Limit Theorem for finitely susceptible
FKG systems}
\author{Oliver Johnson}
\date{\today}

\maketitle

\begin{abstract} 
We adapt arguments concerning entropy-theoretic convergence from the 
independent case to the case of
FKG random variables. FKG systems are chosen since their dependence structure
is controlled through covariance alone, though in the sequel we use many of
the same arguments for weakly dependent random variables. As in previous work
of Barron and Johnson, we consider random variables perturbed by small normals,
since the FKG property gives us control of the resulting densities. We need
to impose a finite susceptibility condition -- that is, the covariance 
between
one random variable and the sum of all the random variables should remain
finite.
\end{abstract}

\section{Introduction and notation}
Gnedenko and Korolev \cite{gnedenko2} discuss the relationship 
between probabilistic limit theorems and the increase of entropy, saying that
\begin{quote}
The formal coincidence of definitions of entropies in physics and in 
information theory gives rise to the question, whether analogs of the second
law of thermodynamics exist in probability theory.
\end{quote}
It is indeed striking that whilst the principle of increase of physical
entropy is taken for granted, the increase of information theoretic
entropy is much less well understood.
The Gaussian is both the distribution of maximum 
entropy (under a variance constraint) and the limit distribution of 
convolutions in the Central Limit regime (which preserves variance). There
is clear physical interest in asking questions such as whether entropy always
increases on convolution, and whether it tends to this maximum.

By showing that the entropy tends to its maximum, we prove the Central
Limit Theorem in a stronger sense than classical methods achieve.
This form of convergence implies the classical weak convergence proved
in Theorem 2 of Newman \cite{newman80}.
Lemma 5 of Takano \cite{takano2} and Theorem 3.1 of Carlen and Soffer 
\cite{carlen} also only prove weak convergence, though under different
conditions. Furthermore by understanding
how score functions become more linear on convolution, we gain an
insight into the workings of the limit theorem, and why the Gaussian 
is the limiting distribution. We are able to gain some insight into
the relationship between maximum entropy distributions and limit theorems
in this way, and see the Central Limit Theorem in a new light.

Gnedenko  and   Korolev  propose   a  programme  to   investigate  the
relationship  between  results  like  the Central  Limit  Theorem  and
maximum entropy  distributions.  This programme has  been developed by
Brown \cite{brown},  Barron \cite{barron}, Johnson  \cite{johnson} and
Barron  and Johnson \cite{johnson5},  who have  used information-theoretic
techniques to prove the Central  Limit Theorem. These papers only deal
with  the  case  of  independent random  variables,  \cite{brown}  and
\cite{barron} in  the case  of identically distributed  variables, and
\cite{johnson}   and  \cite{johnson5}   for   non-identical  variables
satisfying a Lindeberg-like condition.

This paper extends these results and develops new techniques to consider the
case of dependent random variables satisfying the FKG inequalities.
The fact that proofs of 
entropy-theoretic convergence have only previously existed for
independent variables is unfortunate, particularly given
the natural physical interest in dependent systems. In extending Barron's
work, we have shown the link between physical and information-theoretic 
entropies holds in more generality than Gnedenko and Korolev discussed.
\begin{definition} A set of random variables $\{ X_1, X_2, \ldots X_m \}$ is
said
to be FKG if $\cov (F(X_1, X_2, \ldots, X_m), G(X_1, X_2, \ldots X_m)) \geq 
0$ for all increasing functions $F,G$. \end{definition}

FKG (Fortuin-Kastelyn-Ginibre) inequalities hold for many physical
models with `positive 
correlation', as discussed by Newman \cite{newman80}. For example, in the
Ising model with Hamiltonian $H = \sum_{j,k} J(j-k) X_j X_k - h \sum X_j$,
the FKG inequalities hold if $J(r) \geq 0$ for all $r$. Further, FKG 
inequalities hold for percolation models, where $X_j = I(\mbox{vertex $j$ 
is in an infinite cluster})$ and Yukawa models of Quantum Field Theory.

To obtain convergence in relative entropy we use de Bruijn's
identity, which relates the relative entropy to Fisher information of 
perturbed random variables, which have densities we can control.
\begin{definition}
For a random variable $U$ with variance $\sigma^2$ and smooth density $f$, 
we consider the score
function $\rho(u) = f'(u)/f(u)$, the Fisher information $J(U) = \ep \rho^2(U)$,
and the standardised Fisher Information $\fish(U) = 
\sigma^2 J(U) -1 = \ep (\sigma \rho(U) +U/\sigma)^2 \geq 0$.
\end{definition}
\begin{lemma}[de Bruijn]
For $U$ with mean 0, variance 1 and density 
$f$, the relative entropy distance $D(f \| \phi)$ from a standard normal 
can be expressed 
in terms of the Fisher information of $U$ perturbed by normals 
$\rat{Z} \sim N(0, \tau)$:
$$ D(f \| \phi) = \frac{1}{2} \int_0^{\infty} \left( J( U + \rat{Z}) - 
\frac{1}{1+\tau} \right) d\tau 
= \frac{1}{2} \int_0^{\infty} \frac{\fish( U + \rat{Z})}{1+\tau} d\tau.$$
\end{lemma}
Lemma 3 of Newman \cite{newman80} shows that for 
$(S,T)$ FKG, we can control \\ $\cov(f(S), g(T))$.
In our case this is useful because for 
$\rat{Z}_S$, $\rat{Z}_T$ are normal $N(0,\tau)$, independent
of $S,T$ and each other, this means we can control the densities
$p_{X,Y} - p_X p_Y$, where $(X,Y) = (S+ \rat{Z}_S, T+\rat{Z}_T)$.
See Lemma \ref{lem:fkgdensbound} for a discussion of these methods.

First we establish conditions under which convergence $\fish(U) \rightarrow 0$
holds, which implies more conventional forms of convergence: 
\begin{lemma}[Shimizu \cite{shimizu}] \label{lem:shim}
If $U$ has variance $\sigma^2$, density $f$ and distribution function
$F$ then denoting the density and distribution function of a $N(0, \sigma^2)$
by $\phi$ and $\Phi$ respectively:
\begin{eqnarray*}
\sup_u |F(u) - \Phi(u)| \leq 
\int | f(u) - \phi(u) | du & \leq & 4 \sqrt{3} \sqrt{\fish(U)} \\
\sup_u | f(u) - \phi(u) | & \leq & \left( 1 + \sqrt{\frac{6}{\pi}} 
 \right) \sqrt{\fish(U)}
\end{eqnarray*}
\end{lemma}
Indeed weak convergence implies that
$ \lim_{n} \ep h(S_n) = \ep h(Z),$
for all bounded uniformly continuous functions $h$. Convergence in relative
entropy extends this to the class of measurable functions bounded by some
multiple of $x^2 +1$ (see Barron \cite{barron} for further details).
\begin{definition}
Consider a stationary $d$-dimensional system of random 
variables $X_{\vc{u}}$ (where the index $\vc{u} \in {\Bbb Z}^d$), 
with mean zero and finite variance. 
For a particular vector $\vc{x} = (u_1, u_2, \ldots, u_d)$,
we define the box 
$$B_{\vc{u}} = \{ \vc{y}: 0 \leq y_i \leq u_i \mbox{ for all $i$} \}.$$
with volume $|\vc{u}| = |B_{\vc{u}}| = \prod_i u_i$.
We can define $v(\vc{x}) = \var(\sum_{\vc{u} \in B_{\vc{x}}} 
X_{\vc{u}})$ and $U_{\vc{x}} = (\sum_{\vc{u} \in B_{\vc{x}}} X_{\vc{u}} )
/\sqrt{|\vc{x}|}$.
Define
perturbed random variables $\rat{Y}_{\vc{u}} = X_{\vc{u}} + \rat{Z}_{\vc{u}}$,
for $\rat{Z}_\vc{u}$
a sequence of $N(0,\tau)$ independent of $X_{\vc{u}}$ and each other. 
We introduce $\rat{V}_{\vc{x}} = (\sum_{\vc{u} \in B_{\vc{x}}} 
\rat{Y}_{\vc{u}} )/\sqrt{|\vc{x}|} \sim U_{\vc{x}} + \rat{Z}$.
\end{definition}
\begin{condition}[Finite Susceptibility] \label{cond:finsus}
$$ v= \sum_{\vc{u}} \cov(X_\vc{0},X_\vc{u}) < \infty.$$ \end{condition}
\begin{definition} For function $\psi$, define the class of random variables
$X$ with variance $\sigma^2$ such that:
$$ {{\cal{C}}}_\psi = \{ X: \ep X^2 I(|X| \geq R\sigma ) \leq \sigma^2 \psi(R)
\mbox{ for all $R$} \} .$$ \end{definition}
\begin{condition}[Uniform Integrability] \label{cond:unifint}
There exists $\psi$ such that $\rat{V}_{\vc{u}} \in {{\cal{C}}}_\psi$ 
for all $\vc{u}$.
\end{condition}
\begin{theorem} \label{thm:fkgconv}
Consider a stationary collection of
mean zero, finite variance
random variables $X_{\vc{u}}$ obeying the FKG inequalities and finite 
susceptibility (Condition \ref{cond:finsus}). Then 
$$\lim_{n \tends} \left( \sup_{\vc{x}: |\vc{x}| = n} 
\fish(\rat{V}_{\vc{x}}) \right) = 0,$$
if and only if Condition \ref{cond:unifint} (Uniform Integrability) holds.
\end{theorem}
Condition \ref{cond:unifint}
(for stationary FKG variables with finite variance) is actually 
implied by Condition \ref{cond:finsus}. This follows by Newman's proof
\cite{newman80} 
that these conditions imply the Central Limit Theorem, since if $F_n(x)$ is
the distribution function of $\rat{V}_n$, then $F_n(x) \rightarrow \Phi(x)$,
so $\int z^2 I(|z| > N) dF_n(z)  =  1 -\int z^2 dF_n(z) I(|z| \leq N) dz
 \rightarrow 1 -\int z^2 I(|z| \leq N) d\Phi(z) =
\int z^2 I(|z| > N) d\Phi(z)$. 
Carlen and Soffer claim on Page 369
that if $\cov(X_{\vc{0}}, X_{\vc{i}})$ decays at
a rate of $|\vc{i}|^{-t}$, where $t > 2d$, $d$ the dimension of the lattice,
then Condition \ref{cond:unifint} will hold. This roughly 
corresponds to requiring
that $\sum_{\vc{i}} \cov(  X_{\vc{0}}, X_{\vc{i}})^{1/2} < \infty$.

{\bf Note}: 
we do not need to assume that the $X_i$ themselves have densities -- 
even if not, by the following Lemma we obtain weak convergence of the
normalized sums of the original variables.
\begin{definition} Define
$ \kappa(n, \tau) = \sup_{|\vc{u}| \geq n} \fish(\rat{V}_{\vc{u}}).$
\end{definition}
\begin{condition} \label{cond:inte}
For some $n$, $\int \kappa(n, \tau)/(1+\tau) d\tau$ is finite.
\end{condition}
\begin{theorem} 
Consider a stationary collection of
mean zero, finite variance
random variables $X_{\vc{u}}$ with densities, obeying the FKG inequalities and 
Conditions \ref{cond:finsus} and \ref{cond:inte}. Then if $g_{\vc{u}}$ is the
density of $\rat{V}_{\vc{u}}$, then:
$$ D(g_{\vc{u}} \| \phi) \rightarrow 0,$$
if and only if Condition \ref{cond:unifint} (Uniform Integrability) holds. 
\end{theorem}
\begin{proof} Via monotone convergence: $\kappa(n, \tau)$ converges
monotonically to $0$ in $n$, and hence $\int \kappa(n, \tau) d\tau$ converges
to zero. \end{proof}
Newman claims that if instead of scaling by $|\vc{x}|$, we scale by 
$v(\vc{x})$, Condition \ref{cond:finsus} can be relaxed to
Condition \ref{cond:unifint} and Condition \ref{cond:slowvar}. 
\begin{condition} \label{cond:slowvar}
If $K(R) = \sum_{|j| \leq R} \cov(X_0, X_j)$, then
$K(R)$ is slowly varying (that is, for any $\lambda$, $\lim_{R \tends}
K(\lambda R)/K(R) =1$). \end{condition}
He remarks that Condition \ref{cond:unifint} can be checked if for example
$\ep (\rat{V}_n)^4 \leq 3 \left( \ep (\rat{V}_n)^2 \right)^2$, which itself
holds in many cases as a consequence of results such as the Lebowitz 
inequality 
\cite{lebowitz} or the GHS inequality \cite{newman75}.

Takano \cite{takano}, \cite{takano2} deals with the behaviour
of entropy, under a $\delta_4$-mixing condition which seems hard
to check in most useful cases, since it is defined by a ratio of densities.
Further, Takano only proves
convergence of in relative entropy of the `rooms' (in Bernstein's terminology),
equivalent to weak convergence of the original variables. 
Our conclusion holds in the  stronger sense of relative entropy convergence
of the full sequence.
Another paper to use entropy-theoretic methods in the dependent case
is by Carlen and Soffer \cite{carlen}. They use a variety of conditions
which are different to ours, but again only prove weak convergence for
dependent variables.
\section{Fisher Information and convolution}
In the independent case, Fisher information is a sub-additive quantity on
convolution. In the dependent case, we prove that Fisher information is
`almost sub-additive' -- the interest comes in trying to bound the error
term. Takano \cite{takano}, \cite{takano2} produces bounds 
which depend on his $\delta_4$ mixing coefficient, which is hard to 
understand, and hard to check since it depends on ratios of densities. Our 
calculations provide weaker, and more standard conditions under which the
CLT will hold in the sense of convergence of Fisher Information.

\begin{definition}
For random variables $X$, $Y$ with score functions $\rho_X$, $\rho_Y$, for
any $\beta$, we define $\widetilde{\rho}$ for the score function of
$\sqrt{\beta} X + \sqrt{1-\beta} Y$ and then:
$$ \Delta(X, Y, \beta) =
\ep \left( \sqrt{\beta} \rho_X(X) + \sqrt{1-\beta}
\rho_Y(Y) - \widetilde{\rho}\left( \sqrt{\beta} X + \sqrt{1-\beta} Y \right) 
\right)^2 \geq 0.$$
\end{definition}

The principal theorem of this section is:

\begin{theorem} \label{thm:fkgfishsub}
Let $S$ and $T$ be FKG random variables, with mean zero and variance $
\leq K$. Consider $\rat{Z}_S$ and $\rat{Z}_T$, distributed as $N(0,\tau)$,
indepdendent of $S,T$ and of each other.
Define $X = S + \rat{Z}_S$ and 
$Y = T + \rat{Z}_T$, with score functions $\rho_X$ and $\rho_Y$. There
exists a constant $C= C(K, \tau, \epsilon)$ such that for any $\beta$:
$$ \beta J(X) + (1-\beta) J(Y) - J \left(\sqrt{\beta} X + 
\sqrt{1-\beta} Y \right) + C \cov(S,T)^{1/3-\epsilon}
\geq  \Delta(X, Y, \beta).$$
If $S,T$ have bounded $(2+\delta)$th moment, we can replace $1/3$ by
$(2+\delta)/(6+\delta)$.
\end{theorem}
\begin{proof} The proof of the first result 
requires some involved analysis, and is deferred 
to Section \ref{sec:mainproof}. \end{proof}

Next, we need lower bounds on the term
$ \Delta(X, Y, \beta)$.
As discussed in Barron and Johnson \cite{johnson5}, in the case of independent
variables, such terms are equal to zero exactly when all the functions 
concerned are linear. In general, if such an expression is small, then the
derivatives of $\rho_X$ and $\rho_Y$ are close to constant, so long as we
have uniform control over the tails of $X$ and $Y$. 

\begin{proposition} \label{prop:deltadom}
For any $\psi$, there exists a function $\nu = \nu_{\psi}$,
with $\nu(\epsilon) \rightarrow 0$ as $\epsilon \rightarrow 0$, such that
if $X,Y$ lie in ${\cal{C}}_\psi$, then
$$ \beta(1-\beta) \fish(X) \leq \nu \left( \Delta(X, Y, \beta) \right).$$
\end{proposition}
\begin{proof} 
Define a semi-norm $\| \; \|_{\Theta}$ on functions via:
$$ \| f \|_{\Theta}^2 = \inf_{a,b} \; \; 
\ep \left( f(Z_{\tau/2}) - aZ_{\tau/2} -b
\right)^2,$$
where $Z_{\tau/2}$ is $N(0,\tau/2)$. 

Using Lemma 3.1 of Johnson \cite{johnson3},
for $K >0$, there exists a constant $\xi_{K} > 0$ such that for 
any dependent random variables $(S,T)$ with variances $\leq
K$ then the sum $ (X,Y) = (S+ \rat{Z}_S, T+\rat{Z}_T)$ has joint
density $\rat{p}(x,y)$  
bounded below by $\xi_{K} \phi_{\tau/2}(x) \phi_{\tau/2}(y)$.

Hence writing $h(x,y) = \sqrt{\beta} \rho_X(x) + \sqrt{1-\beta}
\rho_Y(y) - \widetilde{\rho} \left( \sqrt{\beta} x + \sqrt{1-\beta} y 
\right)$, then:
\begin{eqnarray*}
\Delta(X,Y,\beta)
& = & \int \rat{p}(x,y) h(x,y)^2 dx dy 
\geq  \xi_{K} \int \phi_{\tau/2}(x) \phi_{\tau/2}(y) h(x,y)^2 dx dy \\
& \geq & \frac{\beta(1-\beta) \xi_{K}}{2} \left(
\| \rho_X  \|_{\Theta}^2 + \| \rho_Y \|_{\Theta}^2 \right),
\end{eqnarray*}
by Proposition 3.2 of \cite{johnson}.
The crucial result of Johnson \cite{johnson} implies
for a fixed $\psi$, if the sequence $X_n \in {{\cal{C}}}_\psi$ have 
score functions $ f_n$, then $\| f_n \|_{\Theta} \rightarrow 0$
implies that $\fish(X_n) \rightarrow 0$. \end{proof}

Define
$$J(n) = \sup \{ \fish( \rat{V}_{\vc{x}}) : |\vc{x}| = n \}.$$
Note that in the 1-dimensional case, there is only one set of this 
form, $\{ 1, 2, \ldots n \}$. 

\begin{corollary} \label{cor:fkgmain}
If Condition \ref{cond:unifint} holds for some $\psi$ then
there exists $d(m) \rightarrow 0$ as $m \rightarrow \infty$ such that
for $m \geq n$:
$$ J(n+m) \leq \frac{m}{m+n} J(m) + \frac{n}{m+n} J(n) + d(m)
- \nu^{-1}_\psi \left( \frac{J(m) mn}{(m+n)^2} \right).$$
\end{corollary}
\begin{proof}
For any $\vc{x}$, we can decompose the box into smaller distinct ones:
$ B_{\vc{x}} = B_{\vc{y}} \bigcup \widetilde{B}_{\vc{z}},$
where $B_{\vc{y}} \cap \widetilde{B}_{\vc{z}} = \emptyset$, and 
$\vc{x} = \vc{y} = \vc{z}$ for all but the $j$th coordinate, so that:
$$ B_{\vc{z}} =  \{ \vc{u}: 0 \leq u_i \leq x_i \mbox{ for all $i \neq j$ 
and } y_j+1 \leq u_j \leq y_j + z_j = x_j \}.$$ 
This corresponds to splitting the box into two smaller ones by making a cut
parallel to the $j$th face. We write $\widetilde{U}_{\vc{z}} = 
(\sum_{\vc{u} \in \widetilde{B}_{\vc{z}}} X_{\vc{u}})/
\sqrt{|\vc{z}|}$, and $\vt_{\vc{z}} = 
(\sum_{\vc{u} \in \widetilde{B}_{\vc{z}}} \rat{Y}_{\vc{u}})/\sqrt{|\vc{z}|}$

Taking $\beta = m/(m+n) = |\vc{y}|/|\vc{x}|$, and by substituting in 
Theorem \ref{thm:fkgfishsub}, since $J(m) \leq 1/\tau$, we obtain
\begin{eqnarray*}
\fish(\rat{V}_{\vc{x}}) & \leq & \frac{m}{m+n} 
\fish(\rat{V}_{\vc{y}}) + \frac{n}{n+m} 
\fish(\rat{V}_{\vc{z}}) \\
& & + C'\cov(\rat{V}_{\vc{y}}, \vt_{\vc{z}})^{1/3-\epsilon} - 
\Delta \left( \rat{V}_{\vc{y}}, \vt_{\vc{z}}, \frac{m}{m+n} \right).
\end{eqnarray*}
Define $$c(m,n) = \sup \left\{
\cov( U_{\vc{y}}, \widetilde{U}_{\vc{z}}) : |\vc{y}|=m, |\vc{z}|=n
\right\}.$$
Under the finite susceptibility condition  (Condition \ref{cond:finsus}),
Lemma 4 of Newman \cite{newman80} shows that this quantity is bounded above
in a suitable way, since $\var(U_{\vc{u}} - U_{\vc{v}}) \rightarrow 0$ if
$|\vc{u}|/|\vc{v}| \rightarrow 0$.
%
%there exists a sequence $d(m)$, where $d(m) \rightarrow 0$ as $m \rightarrow
%\infty$ such that $C c(m,n)^{1/3-\epsilon} \leq d(m)$ when $m \leq n$.
%For any $\vc{y}, \vc{z}$, we obtain:
%\begin{eqnarray*}
%\cov(\rat{V}_\vc{y}, \vt_\vc{z}) &= & 
%\frac{1}{\sqrt{|\vc{x}| |\vc{y}|}} 
%\left( \sum_{\vc{u} \in B_{\vc{y}}} \sum_{\vc{v} \in \widetilde{B}_{\vc{z}}}
%C(|\vc{v} - \vc{u}|) \\
%& \leq & \sqrt{\frac{|\vc{x}|}{ |\vc{y}|}} 
%\left( \sum_{|\vc{u}| \leq r} C(\vc{u}) + \sum_{|\vc{u}| \geq r} C(\vc{u}) \\
%& \leq & v/\sqrt{m} + \sum_{\sqrt{m}+1}^{\infty} C(i) = d(m),
%\end{eqnarray*}
%where $d(m) \rightarrow 0$ as $m \tends$.
\end{proof}

We are able to complete the proof of the CLT, under FKG conditions.

\begin{proof}{\bf of Theorem \ref{thm:fkgconv}}
We first establish convergence along the `powers of 2 subsequence'.
Condition \ref{cond:unifint} implies that $\rat{V}_n 
\in {\cal{C}}_{\psi}$ for some
$\psi$ and hence that $J(\rat{V}_n) \leq \nu_{\psi}(\Delta(\rat{V}_n, 
\vt_n, 1/2))$. We can write $D(k) = \nu^{-1}(\fish(\rat{V}_{2^k}/4))$.
By Corollary \ref{cor:fkgmain}, we know that:
$$ J(2^{k+1}) \leq J(2^k) + d(k) - D(k),$$ 
where $d(k) \rightarrow 0$.

We use an argument structured like Linnik's proof \cite{linnik}.
Given $\epsilon$, we can find $K$ such that $d(k) \leq \epsilon/2$, for all
$k \geq K$. Now either:
\begin{enumerate}
\item{For all $k \geq K$, $2d(k) \leq D(k)$, so 
$ J(2^k) - J(2^{k+1}) \geq D(k)/2,$ 
and summing the telescoping sum, we deduce that
$ \sum_k D(k)$ is finite, and hence there exists $L$ such 
that $D(L) \leq \epsilon$.}
\item{Otherwise for some $L \geq K$, $2d(L) \geq D(L)$, then $D(L) \leq 
\epsilon$.}
\end{enumerate}
Thus, in either case, there exists $L$ such that 
$D(L) \leq \epsilon$, and hence by Proposition 
\ref{prop:deltadom}, $J(2^L) \leq 4\nu(\epsilon).$

Now, for any $k \geq L$, either $J(2^{k+1}) \leq J(2^k)$, or
$D(k) \leq d(k) \leq \epsilon$. In the 
second case, $J(2^k) \leq 4\nu(\epsilon)$, so that
$J(2^{k+1}) \leq 4\nu(\epsilon) + \epsilon$. In either case, we prove
by induction that for all $k \geq L$, that
$J(2^{k+1}) \leq 4\nu(\epsilon) + \epsilon$. 

Now, we can `fill in the gaps' to gain control of the whole sequence, adapting
the proof of the standard sub-additive inequality, using the methods 
described in Appendix 2 of \cite{grimmett2}.
%; given
%$\epsilon$, there exists $T$ such that $J(\rat{V}_{2^t}) \leq \epsilon$
%for $t \geq T$. Thus for $k= 2^{t_1} + \ldots 2^{t_r}$, where $T < t_1
%< \ldots t_r$, we can show by induction that
%\begin{eqnarray*}
% J(\rat{V}_k) & \leq & \epsilon + \sum_{s=2}^{r} \delta(2^{t_1} +
%\ldots 2^{t_{s-1}}, 2^{t_{s}}) \frac{2^{t_1} + \ldots 2^{t_{s}}}{k}
%\leq \epsilon + d(2^T) \sum_{s=2}^{r} \frac{2^{t_1} + \ldots 2^{t_{s}}}{k},
%\end{eqnarray*}
%giving the control we require.
\end{proof}
\section{Proof of sub-additive relation} \label{sec:mainproof}
This is the key part of the argument, proving the  bounds at the heart of
the limit theorems. However, although the analysis is somewhat involved,
it is not too technically difficult.

We introduce notation where it will be clear whether densities or score
functions are associated with joint or marginal distributions, by their
number of arguments: $\rho_X(x)$ will be the score function of $X$, 
and $p'_X(x)$ the derivative of its density. For joint densities
$p_{X,Y}(x,y)$, $\der{p}{1}{X,Y}(x,y)$ will be the derivative of the 
density with respect to the first argument and $\der{\rho}{1}{X,Y}(x,y)
= \der{p}{1}{X,Y}(x,y)/p_{X,Y}(x,y)$, and so on.

Note that a similar equation to the independent case tells us about the 
behaviour of Fisher Information of sums: 
\begin{lemma}
If $X$, $Y$ are random variables, 
with joint density $p(x,y)$, and score functions $\der{\rho}{1}{X,Y}$ and 
$\der{\rho}{2}{X,Y}$ then $X+Y$ has score function $\widetilde{\rho}$
given by
$$ \widetilde{\rho}(z) = 
    \ep \left[ \left. \der{\rho}{1}{X,Y}(X,Y) \right| X+Y=z \right] =
    \ep \left[ \left. \der{\rho}{2}{X,Y}(X,Y) \right| X+Y=z \right].$$
\end{lemma}
\begin{proof}
Since $X+Y$ has density $p_{X+Y}$ given by
$ p_{X+Y}(z) = \int p(z-y,y) dy$, then: 
$$p_{X+Y}'(z) = \int \frac{\partial p}{\partial z}(z-y,y) dy.$$
Hence dividing, we obtain that:
$$ \widetilde{\rho}(z) = \frac{p_{X+Y}'(z)}{p_{X+Y}(z)}   = 
\int \der{\rho}{1}{X,Y}(z-y,y) \frac{p(z-y,y)}{p_{X+Y}(z)} dy,$$
as claimed. \end{proof}

For given $a,b$, define the function $M(x,y) = M_{a,b}(x,y)$ by:
$$ M(x,y) = a \left( \der{\rho}{1}{X,Y}(x,y)- \rho_X(x) \right) +  
b \left( \der{\rho}{2}{X,Y}(x,y)- \rho_Y(y) \right),$$ which is zero 
if $X$ and $Y$ are independent. We will show that if $\cov(X,Y)$ is small, 
then $M$ is close to zero.
\begin{proposition} \label{prop:fishdecomp}
If $X,Y$ are random variables, with score functions $\rho_X,\rho_Y$, 
and if the sum $\sqrt{\beta} X + \sqrt{1-\beta} Y$ has score function 
$\widetilde{\rho}$ then 
\begin{eqnarray*}
\lefteqn{ \beta J(X) + (1-\beta) J(Y) - J \left(\sqrt{\beta} X + 
\sqrt{1-\beta} Y \right)} \\
& &  + 2 \sqrt{\beta(1-\beta)} \ep \rho_X(X) \rho_Y(Y)
+ \ep M_{\sqrt{\beta},\sqrt{1-\beta}}(X,Y) 
\widetilde{\rho}(X+Y)  \\
& = & \ep \left( \sqrt{\beta} \rho_X(X) + \sqrt{1-\beta}
\rho_Y(Y) - \widetilde{\rho} \left(\sqrt{\beta} X + \sqrt{1-\beta} Y \right) 
\right)^2  \end{eqnarray*}
\end{proposition}
\begin{proof}
By the two-dimensional version of Stein's equation, for any function $f(x,y)$:
$$ \ep \der{\rho}{1}{X,Y}(X,Y) f(X,Y) = 
- \ep \frac{\partial f}{\partial x}(X,Y).$$
In particular, if $f(x,y) = \widetilde{\rho}(x+y)$:
$$ \ep \rho_X(X) \widetilde{\rho}(X+Y) = 
- \ep \widetilde{\rho}'(X+Y) - \ep (\der{\rho}{1}{X,Y}(X,Y) - \rho_X(X)) 
\widetilde{\rho}(X+Y).$$
Hence, we know that for any $a,b$:
$$ \ep (a \rho_X(X) + b \rho_Y(Y)) \widetilde{\rho}(X+Y) = 
(a+b) J(X+Y) - \ep M_{a,b}(X,Y) \widetilde{\rho}(X+Y).$$
By considering
$ \int p(x,y) \left( a\rho_X(x) + b\rho_Y(y) - (a+b) \widetilde{\rho}(x+y)
\right)^2 dx dy,$ dealing with the cross term with the expression above,
we deduce that:
\begin{eqnarray*} 
\lefteqn{a^2 J(X) + b^2 J(Y) - (a+b)^2 J (X + Y)}  \\
& & + 2 ab \ep \rho_X(X) \rho_Y(Y)
+ 2 (a+b) \ep M_{a,b}(X,Y) 
\widetilde{\rho}(X + Y)\\
& = & \ep \left( a\rho_X(X) + b \rho_Y(Y) - (a+b)
 \widetilde{\rho}(X + Y) \right)^2 \geq 0.
\end{eqnarray*}
As in the independent case, 
we can rescale, and consider $X' = \sqrt{\beta} X$, $Y' = \sqrt{1-\beta}
Y$, and take $a = \beta, b = 1-\beta$.  
Note that $\sqrt{\beta} \rho_{X'}(u) = \rho_X(u/\sqrt{\beta})$,
$\sqrt{1-\beta} \rho_{Y'}(v) = \rho_Y(v/\sqrt{1-\beta})$. 
\end{proof}
We will show that the two terms on the second line of Proposition 
\ref{prop:fishdecomp} can  be controlled when $(X,Y) = (S+ \rat{Z}_S,
T+ \rat{Z}_T)$, by controlling $\cov(S,T)$ alone. We need control of the 
score functions of perturbed variables. We obtain this in two regions,
firstly in Lemma \ref{lem:scoretail} over the tail, and then in Lemma
\ref{lem:scorebody} over the rest of the real line.

We require an extension of Lemma 3 of Barron \cite{barron} applied to 
single and bivariate random variables:
\begin{lemma} \label{lem:scoretail}
For any random variables $S,T$ as before we define $(X,Y) = (S+
\rat{Z}_S, T+\rat{Z}_T)$ and define $\ra{p}{2\tau}_{U,V}$ for the density of
$(U,V) = (S+\ra{Z}{2\tau}_S,T+\ra{Z}{2\tau}_T)$. 
Now there exists a constant $c_{\tau,k} 
= \sqrt{2} (2k/\tau e)^{k/2}$ such that for all $x$:
\begin{eqnarray*} 
\rat{p}_X(x) |\rho_X(x)|^{k} & \leq & c_{\tau,k} \ra{p}{2\tau}_U(x) \\
\rat{p}(x,y) |\der{\rho}{1}{X,Y}(x,y)|^{k} & \leq & c_{\tau,k} 
\ra{p}{2\tau}_{U,V}(x,y) \\
\rat{p}(x,y) |\der{\rho}{2}{X,Y}(x,y)|^{k} & \leq & c_{\tau,k} 
\ra{p}{2\tau}_{U,V}(x,y) 
\end{eqnarray*}
and hence 
$$ \left( \ep | \rho_X(X) |^k \right)^{1/k}
\leq \sqrt{ \frac{2^{1/k} 2k}{\tau e} }.$$ 
\end{lemma}
\begin{proof} 
We adapt Barron's proof, using H\"{o}lder's inequality 
and the bound; $(u/\tau)^{k} \phi_{\tau}(u) \leq c_{\tau,k} \phi_{2\tau}(u)$ 
for all $u$. 
\begin{eqnarray*}
p'_X(x)^{k} & = & \left( \ep \left( \frac{x-S}{\tau} \right)
\phi_{\tau}(x-S) \right)^{k} \\
& \leq & \left( \ep \left( \frac{x-S}{\tau} \right)^{k} 
\phi_{\tau}(x-S) \right)
\left( \ep  \phi_{\tau}(x-S) \right)^{k-1} \\
& \leq & c_{\tau,k} \left( \ep \phi_{2\tau}(x-S) \right) p_X(x)^{k-1}
\end{eqnarray*}
A similar argument gives the other bounds.
\end{proof}
Now, the normal perturbation ensures that the density does not decrease too
fast, and so the modulus of the score function can not grow too fast.
By considering $S$ normal, so that $\rho$ grows linearly with $u$, we know
that the $B^3$ rate of growth is a sharp bound.
\begin{lemma} \label{lem:scorebody}
If $S$ is a random variable with variance $\leq K$, for 
$X = S + \rat{Z}_S$, with score function $\rho$, for $B>
1$, there exists a function $f_1(\tau,K)$ such that:
$$ \int_{-B\sqrt{\tau}}^{B\sqrt{\tau}} \rho(u)^2 du \leq f_1(\tau,K) B^3.$$
\end{lemma}
\begin{proof} Now: 
$p(u) \geq (2\exp(2K/\tau))^{-1} \phi_{\tau/2}(u)$, so that
for $u \in (-B\sqrt{\tau}, B\sqrt{\tau})$,
$(B\sqrt{\tau} p(u))^{-1} \leq 2\sqrt{\pi} \exp(B^2 + 4/\tau)/B
\leq 2\sqrt{\pi} \exp(B^2 + 4/\tau)
$. Hence for any $k \geq 1$, by H\"older's inequality:
\begin{eqnarray*}
\int_{-B\sqrt{\tau}}^{B\sqrt{\tau}} \rho(u)^2 du & \leq & 
\left( \int_{-B\sqrt{\tau}}^{B\sqrt{\tau}} |\rho(u)|^{2k} du 
\right)^{1/k} \left( 2B \sqrt{\tau} \right)^{1-1/k} \\
& \leq & \left( \int_{-B\sqrt{\tau}}^{B\sqrt{\tau}} 
\frac{ p(u) |\rho(u)|^{2k}}{ 2B \sqrt{\tau} \inf_u p(u)} du 
\right)^{1/k} \left( 2B \sqrt{\tau} \right) \\
& \leq &  \left( \frac{8B}{\sqrt{\tau} e} \right) k \left( 2\sqrt{2\pi} 
\exp(B^2 + 2K/\tau) \right)^{1/k}.  
\end{eqnarray*}
Since we have a free choice of $k \geq 1$ to maximise $k \exp(v/k)$,
since here $v \geq 1$, taking $k=v$  means that $k \exp (v/k) \exp(-1) = v
$. Hence we obtain a bound of
$$ \int_{-B\sqrt{\tau}}^{B\sqrt{\tau}} \rho(u)^2 du \leq  
\frac{8B}{\sqrt{\tau}e} \left( B^2 + \frac{2K}{\tau} + \log(2\sqrt{2\pi}) 
\right) 
\leq \frac{8B^3}{\sqrt{\tau}e} \left( 3+ \frac{2K}{\tau} \right).$$
\end{proof}
We continue by considering $L_B = \{ |x| \leq B \sqrt{\tau}, 
|y| \leq B \sqrt{\tau} \}$.
\begin{lemma} \label{lem:fkgdensbound}
For random variables $S,T$, let  
$X = S + \rat{Z}_S$ and $Y = T + \rat{Z}_T$.
If $S,T$ satisfy the FKG inequalities then 
there exists a function $f_2(\tau,K)$ such that for $B \geq 1$:
$$ \ep M_{a,b}(X,Y) \widetilde{\rho}(X+Y) I( (X,Y) \in L_B) 
\leq  f_2(\tau,K) (a+b) B^4 \cov(S,T).$$
\end{lemma}
\begin{proof}
Lemma 3 of Newman \cite{newman80} uses the fact that FKG inequalities  
imply `positive quadrant dependence', originally due to Lehman \cite{lehman}.
That is, defining $H(s,t) = \pr(S \geq s, T \geq t) - \pr(S \geq s)
\pr(T \geq t)$, $S,T$ are positive quadrant dependent iff
$H(s,t) \geq 0$ for all $s,t$, which is a consequence of $S,T$ being FKG.
Since $\cov(S,T) = \int H(s,t) ds dt \geq 0$, then 
$$\cov (f(S), g(T))
= \int f'(s) g'(t) H(s,t) ds dt \leq \| f' \|_{\infty} \| g' \|_{\infty}
\cov(S,T).$$
Since $|\phi_c(u)'| \leq \exp(-1/2)/\sqrt{2 \pi c^2}$, and 
$|(u \phi_c(u)/c)'| \leq (2 \exp(-3/2))/\sqrt{2 \pi c^3}$, 
we deduce that:
\begin{eqnarray*}
|p_{X,Y}(x,y) - p_X(x) p_Y(y) | & \leq & \frac{\cov(S,T)}{2\pi \tau^2 e},\\
|\der{p}{1}{X,Y}(x,y) - p'_X(x) p_Y(y) | &\leq & \frac{\cov(S,T)}{\pi 
\tau^{5/2} e^2}, \\
|\der{p}{2}{X,Y}(x,y) - p_X(x) p'_Y(y) | & \leq & \frac{\cov(S,T)}{\pi 
\tau^{5/2} e^2}. \\ 
\end{eqnarray*}
We can rearrange $M_{a,b}$ to give
\begin{eqnarray*} 
M_{a,b}(x,y) & = & a \left( \frac{\der{p}{1}{X,Y}(x,y) - p_X'(x) p_Y(y)}
{p_{X,Y}(x,y)} \right) +
b \left( \frac{\der{p}{2}{X,Y}(x,y) - p_X(x) p_Y'(y)}{p_{X,Y}(x,y)} \right)
 \nonumber \\
& & + (a \rho_X(x) + b \rho_Y(y)) \left( 
\frac{p_X(x) p_Y(y) - p_{X,Y}(x,y)}{p_{X,Y}(x,y)} \right).
\end{eqnarray*}
and hence writing $c$ for $\cov(S,T)/(2 \pi \tau^{5/2} e^2)$, f
or $(x,y) \in L_B$:
$$ p_{X,Y}(x,y) |M_{a,b}(x,y)| \leq c 
\left( \sqrt{\tau} e (a\rho_X(x) + b \rho_Y(y)) + 2(a+b) \right).$$
By Cauchy-Schwarz:
\begin{eqnarray*} 
\lefteqn{ \int p_{X,Y}(x,y) M_{a,b}(x,y) \widetilde{\rho}(x+y)
I( (x,y) \in L_B) dx dy } \\
& \leq & c \int \left( \sqrt{\tau} e (a \rho_X(x) + b
\rho_Y(y)) + 2(a+b) \right) \widetilde{\rho}(x+y)
I((x,y) \in L_B) dx dy \\
& \leq & c (a+b) \left( \sqrt{2 B^4 \sqrt{\tau} f_1} \sqrt{16 B^4 
\sqrt{\tau} f_1} + \sqrt{4 B^2 \tau} \sqrt{16 B^4 \sqrt{\tau} f_1} \right)
\end{eqnarray*}

This follows firstly since:
$$ \int \rho_X(x)^2 I((x,y) \in L_B) dx dy  \leq  (2B\sqrt{\tau}) 
\int_{-B\sqrt{\tau}}^{B\sqrt{\tau}} \rho_X(x)^2 dx 
\leq (2B\sqrt{\tau}) B^3 f_1(\tau,K).$$
and
\begin{eqnarray*}
\lefteqn{\int \widetilde{\rho}(x+y)^2 I((x,y) \in L_B) dx dy} \\
 & \leq & \int \widetilde{\rho}(x+y)^2 I(|x+y| \leq 2B\sqrt{\tau})  
I(|y| \leq B\sqrt{\tau})  dx dy \\
& \leq & 2B \sqrt{\tau} 
\int_{-2B\sqrt{\tau}}^{2B\sqrt{\tau}} 
\widetilde{\rho}(z)^2 dz \leq  16B^4 \sqrt{\tau} f_1(\tau,K)
\end{eqnarray*}
\end{proof}

\begin{lemma} \label{lem:scoretailoff}
For any random variables $S,T$ with mean zero and variance 
$\leq K$, let  
$X = S + \rat{Z}_S$ and $Y = T + \rat{Z}_T$.
There exists a function $f_3(\tau,K,\epsilon)$ such that:
$$ \ep M_{a,b}(X,Y) \widetilde{\rho}(X+Y) 
( I( (X,Y) \notin L_B) dx dy \leq (a+b) 
\frac{f_3(\tau,K,\epsilon)}{B^{2-\epsilon}}.$$
for $S,T$ with $k$th moment ($k \geq 2$) bounded above, we can achieve a 
rate of decay of $1/B^{k-\epsilon}$.
\end{lemma}
\begin{proof} By Chebyshev
$\pr \left( (S+\ra{Z}{2\tau}_S,  T+\ra{Z}{2\tau}_T ) \notin L_B) \right)
\leq \int \ra{p}{2\tau}(x,y) (x^2 + y^2)/(2B^2 \tau) dx dy \leq (K+2\tau)/
(B^2 \tau)$
so by H\"older-Minkowski for $1/p + 1/q =1$:
\begin{eqnarray*}
\lefteqn{ \ep \der{\rho}{1}{X,Y}(X,Y) \widetilde{\rho}(X+Y) 
I( (X,Y) \notin L_B)}  \\
& \leq &  \left( \ep | \der{\rho}{1}{X,Y}(X,Y)|^p I((X,Y) \notin L_B)
\right)^{1/p}
\left( \ep | \widetilde{\rho}(X+Y) |^q \right)^{1/q} \\
& \leq & c_{\tau,p}^{1/p} c_{\tau,q}^{1/q}
\pr \left( (S+\ra{Z}{2\tau}_S,  T+\ra{Z}{2\tau}_T ) \notin L_B) \right)^{1/p} 
\\
& \leq & \frac{ 2\sqrt{2} \exp(-1)}{\tau} ( 2+ K/\tau) \sqrt{pq} 
\frac{1}{B^{2/p}}
\end{eqnarray*}
By choosing $p$ arbitrarily close to 1, we can obtain a constant term,
as required. The other terms work in a similar way.
\end{proof}
Similarly we bound the remaining product term:
\begin{lemma} \label{lem:prodbound}
For FKG random variables $S,T$ with mean zero and variance 
$\leq K$, let  
$X = S + \rat{Z}_S$ and $Y = T + \rat{Z}_T$.
There exist functions $f_4(\tau,K)$ and $f_5(\tau,K)$ such that
$$\ep \rho_X(X) \rho_Y(Y)
\leq  f_4(\tau,K) B^4 \cov(S,T) + f_5(\tau,K) /B^{2}. $$
\end{lemma}
\begin{proof}
Using part of Lemma \ref{lem:fkgdensbound}, we know that 
$p_{X,Y}(x,y) - p_X(x) p_Y(y) \leq \cov(S,T)/(2 \pi e \tau^2)$. Hence by 
argument similar to those of Lemmas \ref{lem:scoretailoff} and 
\ref{lem:prodbound}, we obtain that:
\begin{eqnarray*}
\ep \rho_X(X) \rho_Y(Y) 
& = & \int \left( p_{X,Y}(x,y) - p_X(x) p_Y(y) \right)
\rho_X(x) \rho_Y(y) dx dy \\
& \leq & \frac{\cov(S,T)}{2 \pi e \tau^2}
\int  |\rho_X(x)| |\rho_Y(y)| I((x,y) \in L_B) dx dy \\
& & + \int p(x,y) |\rho_X(x)| |\rho_Y(y)| I((x,y) \notin L_B) dx dy \\
& & + \int p(x) p(y) |\rho_X(x)| |\rho_Y(y)| I((x,y) \notin L_B) dx dy  \\
& \leq & \frac{\cov(S,T)}{2 \pi e \tau^2} \left(
\int_{-B\sqrt{\tau}}^{B\sqrt{\tau}} |\rho_X(x)|^2 dx \right)^2 \\
& & + 2 \left( \int p_{X,Y}(x,y) |\rho_X(x)|^2 I((x,y) \notin L_B) dx dy
\right). \\
\end{eqnarray*}
as required.
\end{proof}

\begin{proof}{\bf of Theorem \ref{thm:fkgfishsub}}
Combining Lemmas  \ref{lem:fkgdensbound}, \ref{lem:scoretailoff} and
\ref{lem:prodbound} , we obtain for given $K, \tau, \epsilon$  that
there exist constants $C_1, C_2$ such that
$$\ep M_{\sqrt{\beta},\sqrt{1-\beta}} \widetilde{\rho} 
+ \sqrt{\beta(1-\beta)} \ep \rho_X \rho_Y \leq
C_1 \cov(S,T) B^4 +  C_2/B^{2-\epsilon},$$ so choosing $B = 
(K/\cov(S,T))^{1/6} > 1$, we obtain a bound of 
$C \cov(S,T)^{1/3-\epsilon}$.

By Lemma \ref{lem:scoretailoff}, 
note that if $X,Y$ have bounded $k$th moment, then we obtain decay
at the rate
$C_1 \cov(S,T) B^4 +  C_2/B^{k'}$, for any $k' < k$.
Choosing $B = \cov(S,T)^{-1/(k'+4)}$, we obtain a rate of 
$\cov(S,T)^{k'/(k'+4)}$. Hence
for $k \tends$ we can find a rate arbitarily close to $\cov(S,T)$.
\end{proof}

\end{document}